# SCHOTTKY UNIFORMIZATIONS OF AUTOMORPHISMS OF RIEMANN SURFACES

RUBÉN A. HIDALGO

ABSTRACT. It is well known that the collection of uniformizations of a closed Riemann surface $S$ is partially ordered; the lowest ones are the Schottky unformizations, that is, tuples $(\Omega, \Gamma, P : \Omega \to S)$, where $\Gamma$ is a Schottky group with region of discontinuity $\Omega$ and $P : \Omega \to S$ is a regular holomorphic cover map with $\Gamma$ as its deck group.

Let $\tau : S \to S$ be a conformal (respectively, anticonformal) automorphism of $S$ of finite order $n$, and let $(\Omega, \Gamma, P : \Omega \to S)$ be a Schottky uniformization of $S$. Assume that $\tau$ lifts with respect to the previous Schottky uniformization, that is, there exists a Möbius (respectively, extended Möbius) transformation $\kappa$, keeping $\Omega$ invariant, with $P \circ \kappa = \tau \circ P$. The Kleinian (respectively, extended Kleinian) group $K = \langle \Gamma, \kappa \rangle$ contains $\Gamma$ as a finite index normal subgroup and $K/\Gamma \cong \mathbb{Z}_n$. We provide a structural picture of $K$ in terms of the Klein-Maskit's combination theorems and some basic groups. Some consequences are (i) the determination of the number of topologically different types of such groups (fixed $n$ and the rank of the Schottky normal subgroup) and (ii) for $n$ prime, the number of normal Schottky normal subgroups, up to conjugacy that $K$ has.

As Schottky groups provide those hyperbolic structures of handlebodies whose injectivity radius is bounded away from zero, the structural picture also permits to obtain a geometrical description of finite cyclic groups of homeomorphisms of handlebodies (a 3-dimensional object) from a planar point of view.

## 1. INTRODUCTION

Riemann surfaces, first studied by B. Riemann, are the natural objects on which to do complex analysis in one variable. These are 1-dimensional complex manifolds, that is, they can be though as patching together open pieces of the complex plane with holomorphic transition functions. In this paper we restrict ourselves to the case of closed Riemann surfaces, so they can be described either by (i) smooth complex algebraic projective curves (by the Riemann-Roch Theorem), (ii) principally polarized abelian varieties (by Torelli's Theorem) and (iii) Kleinian groups (by Poincaré-Klein-Koebe's Uniformization Theorem). In this paper we will consider the last type of representations of closed Riemann surfaces.

A *uniformization* of a closed Riemann surface $S$ is a triple $(\Delta, \Gamma, P : \Delta \to S)$, where $\Gamma$ is a Kleinian group, $\Delta$ is a $\Gamma$-invariant connected component of its region of discontinuity and $P : \Delta \to S$ is a regular planar covering with $\Gamma$ as its group of deck transformations. In [23] Maskit provided a description of all uniformizations (regular planar coverings) in terms of a (not unique) collection of pairwise disjoint loops.

The collection of uniformizations of $S$ is partially ordered in the sense that $(\Delta_1, \Gamma_1, P_1 : \Delta_1 \to S)$ is smaller than $(\Delta_2, \Gamma_2, P_2 : \Delta_2 \to S)$ if there is a covering map $Q : \Delta_2 \to \Delta_1$ so that $P_2 = P_1 \circ Q$. Clearly, a uniformization $(\Delta, \Gamma, P : \Delta \to S)$ is a highest one (with respect to the previous partial ordering) if and only if $\Delta$ is simply-connected. The lowest uniformizations [23, 26] are produced when $\Gamma$ is a Schottky group of rank $g$ (in which case $\Delta$ is the region of discontinuity $\Omega$ of $\Gamma$); we call them *Schottky uniformizations* of $S$.

Let us fix a uniformization $(\Delta, \Delta, P : \Delta \to S)$ and a conformal (respectively, anticonformal) automorphism $\tau$ of $S$. We say that $\tau$ *lifts with respect to the above uniformization* if there is an automorphism $\kappa$ of $\Delta$

---

2000 *Mathematics Subject Classification.* Primary 30F10, 30F40.
Partially supported by Projects Fondecyt 1110001 and UTFSM 12.11.01.

so that $P \circ \kappa = \tau \circ P$. Clearly, if $\tau$ lifts, then any power of it also lifts. A group $H$ of (conformal/anticonformal) automorphisms of $S$ *lifts with respect to the above uniformization* if and only if each of its elements lifts. If $S$ has genus at least two, then necessary and sufficient conditions for the lifting property of $H$ were provided in [11] (in Section 3 we recall these conditions for the case of Schottky uniformizations [8]). It was noticed in [6] that, for $S$ of genus at least 2 (so its full group of automorphisms is finite), if $H$ lifts with respect to a uniformization of $S$ which is not a highest one, then the order of $H$ is at most $24(g-1)$ (and if $H$ only consists of conformal automorphisms, then its order is at most $12(g-1)$).

As the region of discontinuity $\Omega$ of a Schottky group is of class $O_{AD}$ (that is, it admits no holomorphic function with finite Dirichlet norm [1, pg 241]), it follows that a conformal (respectively, anticonformal) automorphism of $\Omega$ is the restriction of a Möbius (respectively, extended Möbius) transformation. In this way, if $H$ is finite group of automorphisms of $S$ (only a restriction for genus zero and one) which lifts with respect to a Schottky uniformization of $S$, say $(\Gamma, \Omega, P : \Omega \to S)$, then its lifting provides a (extended) Kleinian group $K$ containing the Schottky group $\Gamma$ as a finite index normal subgroup with $K/\Gamma \cong H$. The finite index property, together with the fact that Schottky groups have no parabolic transformations, asserts that the lifted group $K$ is a (extended) Kleinian group, with $\Omega$ as its region of discontinuity, containing no parabolic transformations. If $H$ only has conformal automorphisms, then Maskit's decomposition theorem of function groups [18, 19] asserts that $K$ can be constructed from a finite collection of finite groups and cyclic subgroups generated by either elliptic or loxodromic transformations by a finite number of applications of the Klein-Maskit combination theorems [20, 21] (similar statements hold in the case that $H$ admits anticonformal automorphisms). As the geometric structure of Schottky groups is simple, one expect to have a simple and concrete decomposition structure for the lifted group $K$ in terms of the algebraic structure of $H$. In the case that $H$ is a group of conformal automorphisms of the maximal order $12(g-1)$, then the structure of $K$ was obtained in [9] and, for $H \cong \mathbb{Z}_2$, it was obtained in [4, 5]. In this paper we extend these previous structure results for the case when $H$ is a cyclic group (see Theorems 1 and 4).

As a direct consequence of our structural decomposition result, we are able to obtain results related to the fixed points and quotients of finite order homeomorphisms of handlebodies (see the final section) [12, 13].

## 2. Main results

Let $S$ be a closed Riemann surface, let $\tau : S \to S$ be a conformal (respectively, anticonformal) automorphism of order $n$ (respectively, $2n$) and let $(\Omega, \Gamma, P : \Omega \to S)$ be a Schottky uniformization of $S$ for which $\tau$ lifts. As already noticed, if $\kappa$ is a lifting of $\tau$, then $K = \langle \kappa, \Gamma \rangle$ is a Kleinian (respectively, extended Kleinian) group with $K/\Gamma$ isomorphic to $\mathbb{Z}_n$ (respectively, $\mathbb{Z}_{2n}$). We say that $K$ is a $\mathbb{Z}_n$-*Schottky group* (respectively, an *extended $\mathbb{Z}_n$-Schottky group*) of rank $g$.

2.1. **Decomposition structure of $\mathbb{Z}_n$-Schottky groups.** The elementary $\mathbb{Z}_n$-Schottky groups are those of rank either 0 or 1. Those of rank 0 are provided by the cyclic groups generated by an elliptic transformation of order $n$ (see Figure 1 with $m = a = 0$ and $b = 1$). The ones of rank 1 are given by (i) cyclic groups generated by a loxodromic transformation (see Figure 1 with $m = b = 0$ and $a = 1$), (ii) groups generated by a loxodromic transformation and an elliptic transformation of order $n_1 \geq 2$ (a divisor of $n$) both of them commuting (see Figure 1 with $m = 1$ and $a = b = 0$) or (iii) the free product of 2 elliptic transformations of order 2, in this case $n = 2$ (see Figure 1 with $a = m = 0$ and $b = 2$). In the general case, we have the following decomposition theorem.

**Theorem 1** (Decomposition theorem for $\mathbb{Z}_n$-Schottky groups). *Every $\mathbb{Z}_n$-Schottky group can be constructed, by use of the Klein-Maskit combination theorem, as free products of the following list of elementary groups (see Figure 1).*

- (T1+) *Cyclic groups generated by loxodromic transformations.*
- (T2) *Cyclic groups generated by elliptic transformations of order a divisor of n.*
- (T4) *Abelian groups generated by a loxodromic transformation and an elliptic transformation of order a divisor of n (in particular, both fixed points of the elliptic are the same as for the loxodromic).*

*Moreover, a Kleinian group constructed as above using "a" groups of type (T1+), "b" groups of type (T2), say cyclic groups of orders $n_1,..., n_a$ (each $n_j > 1$ being a divisor of n), and "m" groups of type (T4), is a $\mathbb{Z}_n$-Schottky group of rank g if conditions (1) and (2) below are satisfied.*

1. $g = n(m + a - 1) + 1 + \sum_{j=1}^{b} \frac{n}{n_j}(n_j - 1)$.
2. *if $a = m = 0$, then $\gcd(n/n_1, ..., n/n_b) = 1$, where gcd stands for "greater common divisor".*

**Remark 2.** The decomposition of $\mathbb{Z}_n$-Schottky groups in Theorem 1 is unique. In fact, let $K$ be a $\mathbb{Z}_n$-Schottky groups (with region of discontinuity $\Omega$), let $O_K^3 = \mathbb{H}^3/K$ be the hyperbolic 3-dimensional orbifold uniformized by $K$ and let $\mathcal{B} \subset O_K^3$ be the locus of cone points (the branch values of the natural quotient map $\mathbb{H}^3 \to O_K^3$). Let us assume that $K$ can be constructed using "a" cyclic groups of type (T1+) "b" cyclic groups of type (T2) and "m" abelian groups of type (T4). Then, the number of connected components of $\mathcal{B}$ which are simple loops is exactly $m$, the number of connected components of $\mathcal{B}$ which are simple arcs is $b$, and $a + m$ is the genus of the conformal boundary $\Omega/K$.

Theorem 1 is proved in Section 4. The main idea will be to construct a suitable $K$-invariant collection of pairwise disjoint simple loops on the region of discontinuity of a given $\mathbb{Z}_n$-Schottky group $K$. That collection of loops will provide the desired decomposition by use of the Klein-Maskit combination theorem (just by using the free product construction as stated in Section 3) of a collection of elementary groups of types (T1+), (T2) and/or (T4).

A direct consequence of Theorem 1 is the following, generalizing the case $p = 2$ obtained in [4].

**Corollary 3.** *Let p be a prime integer. Every Kleinian group containing a Schottky group as an index p normal subgroup can be constructed, by use of the Klein-Maskit combination theorem, as the free product of the following list of elementary groups.*

(i) *cyclic groups generated by either loxodromic transformations or elliptic transformations of order p; and*
(ii) *Abelian groups, isomorphic to $\mathbb{Z} \oplus \mathbb{Z}_p$, generated by a loxodromic transformation and an elliptic transformation of order p (in this case, both share the same fixed points).*

2.2. **Decomposition structure of extended $\mathbb{Z}_n$-Schottky groups.** As for the previous situation, the elementary extended $\mathbb{Z}_n$-Schottky groups are those in rank 0 and 1. These groups are the ones described in the next theorem as (T1), (T2), (T3), (T4), (T5), (T6), (T7) and (T8). These groups are the basic pieces in the construction of all extended $\mathbb{Z}_n$-Schottky groups.

**Theorem 4** (Decomposition theorem for extended $\mathbb{Z}_n$-Schottky groups). *Every extended $\mathbb{Z}_n$-Schottky group can be constructed, by use of the Klein-Maskit combination theorem, as free products of the following list of elementary groups (see Figure 2).*

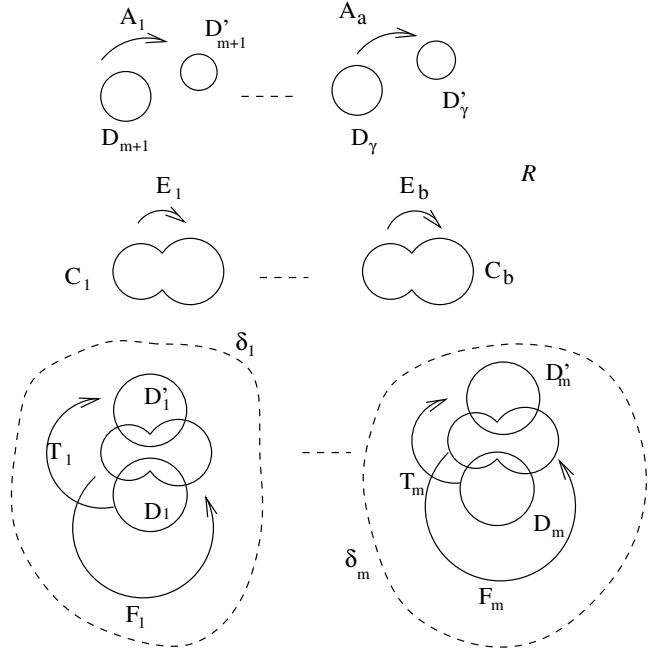

FIGURE 1. Geometric structure of $\mathbb{Z}_n$-Schottky groups ($\gamma = a + m$)

(T1).- *Cyclic groups generated by either a glide-reflection transformation or a loxodromic transformation.*
(T2).- *Cyclic groups generated by elliptic transformation of order a divisor of n.*
(T3).- *Cyclic groups generated by a pseudo-elliptic transformation of order $2d$, where $d$ is a divisor of $n$, but $2d$ is not a divisor of $n$.*
(T4).- *Abelian groups generated by a loxodromic transformation and an elliptic transformation of order a divisor of $n$ (in particular, both fixed points of the elliptic are the same as for the loxodromic).*
(T5).- *Groups generated by a loxodromic transformation $A$ and a pseudo-elliptic transformation $B$, of order a divisor of $2n$ but not of $n$, so that $B^{-1}ABA = I$ (in particular, both fixed points of the loxodromic transformation are permuted by the pseudo-elliptic transformation).*
(T6).- *If $n$ is even, groups generated by a glide-reflection transformation $A$ and an elliptic transformation $B$ of order $2$ so that $BA = AB = I$ (in particular, both fixed points of the elliptic are the same as for the glide-reflection).*
(T7).- *If $n$ is odd, cyclic groups of order $2$ generated by reflections.*
(T8).- *If $n$ is odd, groups generated by the reflection of a circle $\Sigma$ and a discrete group $F$ (of orientation-preserving conformal automorphisms) keeping invariant $\Sigma$ and so that $\Omega_F/F$ (where $\Omega_F$ is the region of discontinuity of $F$) is a connected Riemann orbifold whose conical points have orders divisors of $n$.*

*Moreover, a group constructed using the above basic groups is a extended $\mathbb{Z}_n$-Schottky group if and only if conditions (1) and (2) below are satisfied.*

(1) *In the construction there are groups of either type (T1) or (T3) or (T5) or (T6) or (T7) or (T8); and*

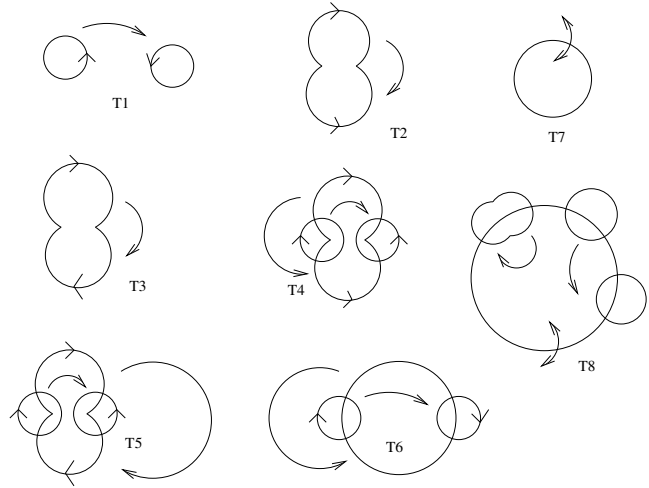

FIGURE 2. Geometric structure of extended $\mathbb{Z}_n$Schottky groups

(2) *If $n \geq 2$ and there are not groups of type (T1) nor (T6) nor groups of type (T8) containing a glide-reflection, then the greater common divisor of all the values of the form $2n/r$, were r runs over all orders of elliptic and pseudo-elliptic transformations used in the basic groups of type (T2), (T3), (T4), (T5), (T7) and (T8), is 1.*

Theorem 4, proved in Section 5, generalizes the decomposition structure of the extended $\mathbb{Z}_1$-Schottky groups (also called *extended Schottky groups*) obtained in [5].

**Corollary 5.** (1) *Every extended Kleinian group, whose index two orientation preserving half subgroup is a Schottky group, can be constructed, by use of the Klein-Maskit combination theorem, as free products of the following list of elementary groups.*
  (i) *cyclic groups generated by reflections;*
 (ii) *cyclic groups generated by imaginary reflections;*
(iii) *cyclic groups generated by glide-reflections;*
(iv) *cyclic groups generated by loxodromic transformations; and*
 (v) *real Schottky groups.*
(2) *Any group K constructed, using a groups of type (i) b groups of type (ii) c groups of type (iii) d of groups of type (iv) and e groups of type (v), is an extended Kleinian group, whose index two orientation preserving half subgroup $\Gamma$ is a Schottky group of rank g, if and only if $a + b + c + e > 0$. Moreover, if the e real Schottky groups are of rank $\gamma_1, ..., \gamma_e \geq 1$, then $\Gamma$ is a Schottky group of rank*

$$g = a + b + 2c + 2d + e - 1 + \sum_{j=1}^{e} \gamma_j.$$

2.3. **Schottky subgroups of index $n$ of $\mathbb{Z}_n$-Schottky groups.** A *geometric automorphism* of a Kleinian group $K$ is an orientation-preserving homeomorphism of the Riemann sphere that conjugates $K$ into itself.

Let $K$ be a $\mathbb{Z}_n$-Schottky group of rank $g$. By the definition, we know the existence of a Schottky group $\Gamma$ of rank $g$ which is a normal subgroup of $K$ with $K/\Gamma \cong \mathbb{Z}_n$, for some $n$. Such a Schottky normal subgroup of $K$ may not be unique. Let us assume $\Gamma_1$ and $\Gamma_2$ are Schottky groups of rank $g$, both normal subgroups of $K$ with $K/\Gamma_j \cong \mathbb{Z}_n$, for $j = 1, 2$. The groups $K$, $\Gamma_1$ and $\Gamma_2$ all have the same region of discontinuity $\Omega$. In this case, for $j = 1, 2$, the closed Riemann surface $S_j = \Omega/\Gamma_j$ admits the cyclic group $H_j = K/\Gamma_j$ as a group of holomorphic automorphisms with $S_j/H_j = O = \Omega/K$. Let us denote by $\pi_j : S_j \to O$ the natural quotient map induced by the action of $H_j$. We say that the Schottky groups $\Gamma_1$ and $\Gamma_2$ are *geometrically equivalent* in $K$ if there exists a geometric automorphism of $K$, say $f : \widehat{\mathbb{C}} \to \widehat{\mathbb{C}}$, so that $f\Gamma_1 f^{-1} = \Gamma_2$. Notice that in this case there are induced orientation-preserving homeomorphisms $\widetilde{f} : S_1 \to S_2$ and $\widehat{f} : O \to O$ so that $\widehat{f} \circ \pi_1 = \pi_2 \circ \widetilde{f}$.

Next result, proved in Section 6, provides the number of geometrically equivalent classes of Schottky subgroups that a $\mathbb{Z}_p$-Schottky group may contain as a normal subgroup of index $p$, for $p$ a prime integer.

**Theorem 6.** *Let $p$ be a prime integer and let $K$ be a $\mathbb{Z}_p$-Schottky group, constructed as in Corollary 3, using "a" loxodromic cyclic groups, "b" elliptic cyclic groups of order $p$, and "m" Abelian groups $\mathbb{Z} \oplus \mathbb{Z}_p$. Then the number of different Schottky normal subgroups of index $p$ that $K$ contains, up to conjugation by a geometrical automorphism of $K$, is exactly*

$$M(b, m) = \binom{b + (p-3)/2}{(p-3)/2} \binom{m + (p-3)/2}{(p-3)/2}$$

Theorem 6 and the fact that the orientation-preserving half of a extended $\mathbb{Z}_p$-Schottky group is a $\mathbb{Z}_p$-Schottky group may be used to extend the above to the case of extended $\mathbb{Z}_p$-Schottky groups. This is not done in this paper.

2.4. **Topologically non-equivalent $\mathbb{Z}_n$-Schottky groups.** The (extended) Kleinian groups are said to be *(weakly) topologically equivalent* if there is a homeomorphism of the Riemann sphere that conjugates one onto the other.

If we fix integer values $g \geq 2$ and $n \geq 2$, then Theorem 1, together Remark 2, permits to obtain the number of topologically non-equivalent $\mathbb{Z}_n$-Schottky groups of rank $g$. In fact, such a theorem asserts that the number of topologically non-equivalent $\mathbb{Z}_n$-Schottky groups of rank $g$ is the same as the number of different tuples of non-negative integers $(m, a, b, n_1, ..., n_a, l_1, ..., l_m)$, where $m, a, b \geq 0$ and $n_j, l_i \geq 2$ are divisors of $n$, satisfying the following conditions:

(1) $$2 \leq n_1 \leq n_2 \leq \cdots \leq n_b \leq n,$$

(2) $$2 \leq l_1 \leq l_2 \leq \cdots \leq l_m \leq n,$$

(3) $$\text{if } m = a = 0, \text{ then } gcd(n/n_1, ..., n/n_b) = 1,$$

and

(4) $$g = n(m + a - 1) + 1 + \sum_{j=1}^{b} \frac{n}{n_j}(n_j - 1).$$

Let $\psi(n)$ be the number of divisors $x \geq 2$ of $n$ and let $2 \leq q_1 < \cdots < q_{\psi(n)} = n$ be all such divisors of $n$. If $m > 0$, then the number $M_m(n)$ of tuples $(l_1, ..., l_m)$, where each $l_j \geq 2$ is a divisor of $n$ satisfying (2), is equal to the number of tuples $(d_1, ..., d_{\psi(n)})$ so that $d_j \geq 0$ and $d_1 + \cdots + d_{\psi(n)} = m$ ($d_j$ denotes the number of repetitions of $l_j$).

Set by $N(n, g; m)$ the number of the possible tuples $(m, a, b, n_1, ..., n_a)$, where $n_j$ are divisors of $n$ satisfying conditions (1), (3) and (4). Equality (4) asserts that $g \geq n(m-1) + 1$, that is, $m \leq (g+n-1)/n$. In this way, all the above together asserts the following fact.

**Corollary 7.** *If $n, g \geq 2$, then there are exactly*

$$N(n, g) = N(n, g; 0) + \sum_{m=1}^{[(g+n-1)/n]} M_m(n) N(n, g; m)$$

*topologically different $\mathbb{Z}_n$-Schottky group of rank $g$.*

2.4.1. *The prime case.* Assume $n = p$ is a prime integer. In this case $n_j = l_k = p$, for $j = 1, ..., b$ and $k = 1, ..., m$. The conditions (1), (2) and (3) are trivially satisfied and $M_m(p) = 1$, for $m > 0$. Now condition (4) is equivalent to

$$b = \frac{g - p(m+a)}{p-1} + 1 \in \{0, 1, ..\},$$

that is, $m + a \in \{0, 1, ..., [(g+p-1)/p]\}$, where $[u]$ denotes the integral part of $u$, and $g - p(m+a) \equiv 0 \mod (p-1)$. In this way, $N(p, g)$ is exactly the number of pairs $(m, a)$ satisfying the following properties:

(i) $m, a \in \{0, 1, ..., [(g+p-1)/p]\}$,
(ii) $m + a \leq [(g+p-1)/p]$,
(iii) $g - p(m+a) \equiv 0 \mod (p-1)$.

If $p = 2$, then every pair $(m, a)$ with $m + a \in \{0, 1, ..., [(g+1)/2]\}$ satisfies (i)-(iii), in particular

$$N(2, g) = \frac{(1 + [\frac{g+1}{2}])(2 + [\frac{g+1}{2}])}{2}.$$

If $p = 3$, then condition (iii) is equivalent to $g - 3(m+a)$ being even. In this case,

$$N(3, g) = \begin{cases} \frac{([\frac{g+2}{3}]+2)^2}{4}, & \text{if } g \text{ is even and } [\frac{g+2}{3}] \text{ is even} \\[6pt] \frac{([\frac{g+2}{3}]+1)^2}{4}, & \text{if } g \text{ is even and } [\frac{g+2}{3}] \text{ is odd} \\[6pt] \frac{[\frac{g+2}{3}]([\frac{g+2}{3}]+2)}{4}, & \text{if } g \text{ is odd and } [\frac{g+2}{3}] \text{ is even} \\[6pt] \frac{([\frac{g+2}{3}]+1)([\frac{g+2}{3}]+3)}{4}, & \text{if } g \text{ is odd and } [\frac{g+2}{3}] \text{ is odd} \end{cases}$$

If $p \geq 5$ is a prime, then $N(p, g)$ and the corresponding triples $(m, a, b)$ can be computed with the help of a computer program. For instance, a short program, written for MATHEMATICA, which can be used to obtain $N(p, g)$ is the following one.

```
n[p_, g_] := Block[{}, k := 0;
Do[ Do[ Do[ If[m + a + b - 1 ≠ 0,
```

```
If[(g + b − 1)/(m + a + b − 1) == p, {k = k + 1,
Print["(", k, ")", "m = ", m, ",", "a = ", a, ",", "b = ", b, ",", "p = ", p]}]],
{b, 0, (g + 1 − 2 ∗ (m + a))}], {m, 0, g}], {a, 0, g}];
Print["N(", p, ",", g, ") = ", k]]
```

One obtain, for instance, that $N(5,5) = 2$, $N(5,10) = 3$, $N(11,10) = 1$, $N(11,100) = 12$, $N(13,157) = 16$. Also, for $p = 11$ and $g = 100$, we get

$$(m, a, b) \in \left\{ \begin{array}{l} (11, 0, 0), (0, 0, 10), (0, 1, 9), (0, 2, 8), (0, 3, 7), (0, 4, 6), \\ (0, 5, 5), (0, 6, 4), (0, 7, 3), (0, 8, 2), (0, 9, 1), (0, 10, 0) \end{array} \right\}.$$

and for $p = 13$ and $g = 157$, we get

$$(m, a, b) \in \left\{ \begin{array}{l} (13, 0, 1), (0, 0, 13), (13, 1, 0), (0, 1, 12), (0, 2, 11), (0, 3, 10), (0, 4, 9), (0, 5, 8), \\ (0, 6, 7), (0, 7, 6), (0, 8, 5), (0, 9, 4), (0, 10, 3), (0, 11, 2), (0, 12, 1), (0, 13, 0) \end{array} \right\}.$$

**Remark 8** (Topologically non-equivalent extended $\mathbb{Z}_n$-Schottky groups). The structure provided by Theorem 4 may also be used to obtain a formula that permits to obtain the number of topologically non-equivalent extended $\mathbb{Z}_n$-Schottky groups. This is a much more complicate formulae than the one provided above for the case of $\mathbb{Z}_n$-Schottky groups. We do not provide it in the general situation, but in Section 7 we work out the case of extended $\mathbb{Z}_2$-Schottky groups.

## 3. Preliminaries

In this section we review some of the definitions (which have not been already stated before), set some notations and recall some technical results we will need in this paper. Generalities on Kleinian and extended Kleinian groups can be found, for instance, in the books [15, 24].

We use the symbol $\Gamma < K$ (respectively, $\Gamma \triangleleft K$) to say that $\Gamma$ is a subgroup (respectively, normal subgroup) of a group $K$. The composition of the maps $f$ and $h$ is as usually denoted by the symbol $f \circ h$, but if we are composing (extended) Möbius transformations $A$ and $B$ we will use the symbol $AB$.

We denote by $\mathbb{M}$ the group of *Möbius transformations* and by $\widehat{\mathbb{M}}$ the group generated by $\mathbb{M}$ and the complex conjugation $J(z) = \bar{z}$. A transformation in $\widehat{\mathbb{M}} - \mathbb{M}$ is called an *extended Möbius transformation*. It is well known that $\mathbb{M}$ is the full group of conformal automorphisms of the Riemann sphere $\widehat{\mathbb{C}}$ and that $\widehat{\mathbb{M}}$ is the full group of conformal and anticonformal automorphisms of it. If $K < \widehat{\mathbb{M}}$, then we set $K^+ = K \cap \mathbb{M}$ and, when $K \neq K^+$, we say that $K^+$ is the *orientation-preserving half* of $K$.

An extended Möbius transformation whose square is an elliptic transformation is called *pseudo-elliptic* (if the square is the identity, then we say that it is a *reflection* if it has fixed points or a *imaginary reflection* otherwise). Similarly, if the square is a loxodromic transformation (in fact a hyperbolic one), then we say that it is a *glide-reflection*. If the square is parabolic, then we say that it is *pseudo-parabolic*.

A *Kleinian group* (respectively, an *extended Kleinian group*) is a discrete subgroup of $\mathbb{M}$ (respectively, a discrete subgroup of $\widehat{\mathbb{M}}$ necessarily containing extended Möbius transformations). The *region of discontinuity* of a (extended) Kleinian group $K$ is the open set (which might be empty) $\Omega \subset \widehat{\mathbb{C}}$ consisting of those points on which $K$ acts discontinuously. The complement closed set $\Lambda = \widehat{\mathbb{C}} - \Omega$ is called the *limit set* of $K$. If $\Lambda$ is finite, then $K$ is called *elementary*; and *non-elementary* otherwise. If $K_1 < K_2 < \widehat{\mathbb{M}}$ and $K_1$ has finite index in $K_2$, then one is discrete if and only if the other is to; in which case both have the same region of discontinuity. In particular, if $K < \widehat{\mathbb{M}}$ contains extended Möbius transformations, then $K$ is an extended Kleinian group if and only if $K^+$ is a Kleinian group. A *function group* (respectively, an *extended function* *group*) is a Kleinian group (respectively, and extended Kleinian group) containing an invariant connected component of its region of discontinuity. Notice that if $K$ is a an extended function group, then $K^+$ is a function group, but the converse is not true in general.

The decomposition of function groups, in the sense of the Klein-Maskit combination theorems is provided in [17, 18, 19] and that for extended function group can be seen in [10] (were the subtle modifications of the arguments for function groups are provided). We next state a simple version of the Klein-Maskit combination theorems which is enough for us in this paper.

**Theorem 9** (Klein-Maskit's combination theorems [20, 21]).
*(1) (Free products) For $j = 1, 2$, let $K_j$ be a (extended) Kleinian group with region of discontinuity $\Omega_j$ and let $\mathcal{F}_j$ be a fundamental domain for $K_j$. Assume that there is a simple closed loop $\Sigma$, contained in the interior of $\mathcal{F}_1 \cap \mathcal{F}_2$, bounding two discs $D_1$ and $D_2$, so that, for $j = 1, 2$, the set $\Sigma \cup D_j \subset \Omega_{3-j}$ is precisely invariant under the identity in $K_{3-j}$. Then $K = \langle K_1, K_2 \rangle$ is a (extended) Kleinian group, with fundamental domain $\mathcal{F}_1 \cap \mathcal{F}_2$, which is the free product of $K_1$ and $K_2$. Every finite order element in $K$ is conjugated in $K$ to a finite order element of either $K_1$ or $K_2$. Moreover, if both $K_1$ and $K_2$ are geometrically finite, then $K$ is so.*
*(2) (HNN-extensions) Let $K_0$ be a (extended) Kleinian group with region of discontinuity $\Omega$ and let $\mathcal{F}$ be a fundamental domain for $K_0$. Assume that there are two pairwise disjoint simple closed loops $\Sigma_1$ and $\Sigma_2$, both of them contained in the interior of $\mathcal{F}_0$, so that $\Sigma_j$ bounds a disc $D_j$ such that $(\Sigma_1 \cup D_1) \cap (\Sigma_2 \cup D_2) = \emptyset$ and with $\Sigma_j \cup D_j \subset \Omega$ precisely invariant under the identity in $K_0$. If $T$ is either a loxodromic transformation or a glide-reflection so that $T(\Sigma_1) = \Sigma_2$ and $T(D_1) \cap D_2 = \emptyset$, then $K = \langle K_0, f \rangle$ is a (extended) Kleinian group, with fundamental domain $\mathcal{F}_1 \cap (D_1 \cup D_2)^c$, which is the HNN-extension of $K_0$ by the cyclic group $\langle T \rangle$. Every finite order element of $K$ is conjugated in $K$ to a finite order element of $K_0$. Moreover, if $K_0$ is geometrically finite, then $K$ is so.*

The *Schottky group of rank* 0 is just the trivial group. A *Schottky group of rank* $g \geq 1$ is a Kleinian group $\Gamma$ generated by loxodromic transformations $A_1, \ldots, A_g$, so that there are $2g$ pairwise disjoint simple loops, $C_1, C'_1, \ldots, C_g, C'_g$, bounding a $2g$-connected domain $\mathcal{D} \subset \widehat{\mathbb{C}}$, where $A_i(C_i) = C'_i$, and $A_i(\mathcal{D}) \cap \mathcal{D} = \emptyset$, for $i = 1, \ldots, g$. The collection of loops $C_1, C'_1, \ldots, C_g$ and $C'_g$ is called a *fundamental set of loops* for $\Gamma$ with respect to the above generators (these groups are constructed, using part (1) in Klein-Maskit's combination theorem, as the free product of cyclic loxodromic groups). The region of discontinuity $\Omega$ of a Schottky group $\Gamma$ of rank $g$ is known to be connected and dense in $\widehat{\mathbb{C}}$ and that $S = \Omega/\Gamma$ is a closed Riemann surface of genus $g$ (the classical retrosection theorem states that, up to conformal isomorphism, every closed Riemann surface is obtained in this way). It is well known that a Schottky group of rank $g$ can be defined as a purely loxodromic Kleinian group of the second kind which is isomorphic to a free of rank $g$ [22].

Let $S$ be a closed Riemann surface, let $(\Delta, \Gamma, P : \Delta \to S)$ be a uniformization of $S$ and let $H < Aut(S)$. In the introduction we have defined the property for $H$ to lift with respect to this uniformization. Note that a lifting $k \in Aut(\Delta)$ of $h \in H$ is not required to be the restriction of a (extended) Möbius transformation. If $K$ is the group generated by all these liftings, then $K < Aut(\Delta)$, $\Gamma \triangleleft K$ and $K/\Gamma \cong H$. As noticed above, in general, $K$ is not a group of (extended) Möbius transformations, but in the case of Schottky uniformizations $K$ turns out to be so. If $\Delta$ is simply-connected, clearly $Aut(S)$ lifts with respect to this uniformization. If $\Delta$ is not simply-connected, it may be that some automorphism of $S$ does not lift to an automorphism of $\Delta$. The following result provides necessary and sufficient conditions for the lifting property to hold in the case of Schottky uniformizations.

**Theorem 10.** [8] *Let $(\Omega, \Gamma, P : \Omega \to S)$ be a Schottky uniformization of the closed Riemann surface $S$ of genus $g \geq 2$. Let $H$ be a group of automorphisms (conformal/anticonformal) of $S$. Then, $H$ lifts with respect* *to the above Schottky uniformization if and only if there is a collection $\mathcal{F}$ of pairwise disjoint simple loops on S such that:*

   (i) *each connected component of $S - \mathcal{F}$ is a planar surface;*
   (ii) *$\mathcal{F}$ is invariant under the action of H; and*
   (iii) *for each $\alpha \in \mathcal{F}$, $P^{-1}(\alpha)$ is a collection of pairwise disjoint simple loops in $\Omega$.*

A collection of loops $\mathcal{F}$, as in Theorem 10, which is minimal (in the sense that by deleting a non-empty sub-collection from $\mathcal{F}$, then one of the above three properties fails) will be called a *fundamental collection of loops* associated to the pair $\{(\Omega, G, P : \Omega \to S), H\}$.

**Remark 11.** We should note that Theorem 10 can be seen a consequence of the Equivariant Loop Theorem [25], whose proof is based on minimal surfaces, that is, surfaces that minimize locally the area. The proof given in [8] only uses arguments proper to (planar) Kleinian groups and the hyperbolic metric.

We will also need the following fact on finite extensions of Schottky groups.

**Theorem 12** ([7]). *Let K be a group of Möbius transformations containing a Schottky group $\Gamma$ as subgroup of finite index. Let h be any elliptic element of K and let x and y be its fixed points. Then either: (i) x and y are both in the region of discontinuity of K or (ii) there is a loxodromic element in $\Gamma$ commuting with h.*

Each Möbius (respectively, extended Möbius) transformation acts (by Poincaré's extension) as an orientation-preserving (respectively, orientation-reversing) isometry of the hyperbolic 3-space $\mathbb{H}^3 = \{(z,t) \in \mathbb{C} \times \mathbb{R} : t > 0\}$ with the hyperbolic metric $ds^2 = (|dz|^2 + dt^2)/t^2$. If $\Gamma$ is a (extended) Kleinian group, then $O_\Gamma^3 = \mathbb{H}^3/\Gamma$ is a 3-dimensional hyperbolic orbifold and the 2-dimensional orbifold $O_\Gamma^2 = \Omega/\Gamma$ is its conformal boundary. In the case that $\Gamma$ is a torsion free Kleinian group, then $O_\Gamma^3$ is a hyperbolic 3-manifold and $O_\Gamma^2$ a Riemann surface. For instance, if $\Gamma$ is a Schottky group of rank $g$, then $O_\Gamma^3$ is homeomorphic to a handlebody of genus $g$ and its conformal boundary $O_\Gamma^2$ is a closed Riemann surface of genus $g$. It is well known that every torsion free Kleinian group $\Gamma$, for which $O^3/\Gamma$ is homeomorphic to the interior of a handlebody of genus $g$ and whose conformal boundary is a closed Riemann surface, is a Schottky group of rank $g$.

### 4. Topological classification of $\mathbb{Z}_n$-Schottky groups

In this section we provide the proof of Theorem 1. We first construct some $\mathbb{Z}_n$-Schottky groups by use of Klein-Maskit's combination theorems with the basic groups of that theorem. Then we show that the constructed groups are all of them.

**4.1. Geometrical $\mathbb{Z}_n$-Schottky groups.** Let us fix some integers $n \geq 2$, $a, b, m \geq 0$, $n_1, ..., n_b, l_1, ..., l_m \in \{2, ..., n\}$ with $n_j$ and $l_i$ divisors of $n$. Set $\gamma = a + m$. Let us consider a collection of $2\gamma + b$ pairwise disjoint simple loops, say $D_1, D'_1,..., D_\gamma, D'_\gamma, C_1,..., C_b$, all of them bounding a common domain $\mathcal{R}$ of connectivity $2\gamma + b$. We also consider $m$ pairwise disjoint simple loops, say $\delta_1, ..., \delta_m \subset \mathcal{R}$, so that, for each $j = 1, ..., m$, $\delta_j$ separates both loops $D_j$ and $D'_j$ from all the others loops (see Figure 1). Let $\Delta_j$ be the closed topological disc bounded by $\delta_j$ which does not contains $D_j$.

Assume there are $\gamma$ loxodromic transformations $A_1,..., A_a, T_1,..., T_m$, and $(b+m)$ elliptic transformations, say $E_1,..., E_b, F_1,..., F_m$, so that

  (i) $E_j$ has order $n_j$, for $j = 1, ..., b$;
  (ii) $F_j$ has order $l_j$, for $j = 1, ..., m$;
  (iii) $A_j(D_{m+j}) = D'_{m+j}$ and $A_j(\mathcal{R}) \cap \mathcal{R} = \emptyset$, for $j = 1, ..., a$;

(iv) $T_j(D_j) = D'_j$ and $T_j(\mathcal{R}) \cap \mathcal{R} = \emptyset$, for $j = 1, ..., m$;
(v) $E_j(C_j) = C_j$ and $E_j(\mathcal{R}) \cap \mathcal{R} = \emptyset$, for $j = 1, ..., b$;
(vi) $F_j$ commutes with $T_j$, for $j = 1, ..., m$; and
(vii) $F_j(\Delta_j) \cap \Delta_j = \emptyset$.

Condition (vii) ensures we may construct a Kleinian group $K$ using the Klein-Maskit combination theorem and the cyclic groups generated by each of the Möbius transformations $A_1,..., A_a, T_1,..., T_m, E_1,..., E_b$, $F_1,...,$ and $F_m$. The same theorem states that

$$K \cong \underbrace{\mathbb{Z} \oplus \mathbb{Z}_{l_1} * \cdots * \mathbb{Z} \oplus \mathbb{Z}_{l_m}}_{m} * \underbrace{\mathbb{Z} * \cdots * \mathbb{Z}}_{a} * \underbrace{\mathbb{Z}_{n_1} * \cdots * \mathbb{Z}_{n_b}}_{b}.$$

We say that the constructed group $K$ above is a *general group of type* $(m, a, b, n_1, ..., n_b, l_1, ..., l_m)$.

Note that, as none of the above cyclic groups contains parabolic transformations, the Klein-Maskit combination theorem asserts that $K$ cannot contain parabolic transformations. It follows that $K$ may only contain elliptic and loxodromic transformations. A fundamental domain for $K$ is drawn in Figure 1. Unfortunately, a Kleinian group as constructed above may fail to be a $\mathbb{Z}_n$-Schottky group. The following provides a necessary and sufficient condition for $K$ to be a $\mathbb{Z}_n$-Schottky group.

**Theorem 13.** *A general group $K$ of type* $(m, a, b, n_1, ..., n_b, l_1, ..., l_m)$ *is a $\mathbb{Z}_n$-Schottky group if and only if (i) or (ii) below holds.*

  (i) $m + a > 0$.
  (ii) $m = a = 0$ and $gcd(n/n_1, ..., n/n_b) = 1$.

*Proof.* A general group $K$ of type $(m, a, b, n_1, ..., n_b, l_1, ..., l_m)$ is a $\mathbb{Z}_n$-Schottky group if and only if there is a surjective homomorphism $\Phi : K \to \mathbb{Z}_n = \langle x : x^n = 1 \rangle$ with torsion free kernel. One direction is clear: if $K$ is a $\mathbb{Z}_n$-Schottky group and $\Gamma \triangleleft K$ is a Schottky group so that $K/\Gamma \cong \mathbb{Z}_n$, then $\Phi : K \to K/\Gamma$ is the canonical projection homomorphism. In the other direction, let us assume we have a surjective homomorphism $\Phi : K \to \mathbb{Z}_n = \langle x : x^n = 1 \rangle$ with torsion free kernel $\Gamma$. As consequence of the Klein-Maskit combination theorem [15], the group $K$ has a connected region of discontinuity, that is, $K$ is a function group. As the group $K$ has no parabolic transformations, the group $\Gamma$ is a torsion free, purely loxodromic function group. That $\Gamma$ is a Schottky group is consequence of the the classification of finitely generated function groups [16].

Now, given a general group $K$ of type $(m, a, b, n_1, ..., n_b, l_1, ..., l_m)$, we need to provide conditions in order to ensure a surjective homomorphism $\Phi$ as above. It is clear that for $\Phi$ to have torsion free kernel it is necessary that $\Phi(F_j)$ and $\Phi(E_j)$ have orders $l_j$ and $n_j$, respectively.

4.2. **Case** $m + a > 0$. In this case set $\Phi : K \to \mathbb{Z}_n = \langle x : x^n = 1 \rangle$ by the rules $\Phi(F_j) = x^{n/l_j}$ $(j = 1, ..., m)$, $\Phi(E_j) = x^{n/n_j}$ $(j = 1, ..., b)$, and
  (i) if $a > 0$, then $\Phi(A_1) = x$, $\Phi(A_j) = \Phi(T_k) = 1$ $(j = 2, .., a, k = 1, ..., m)$,
  (ii) if $a = 0$, then $\Phi(T_1) = x$, $\Phi(T_k) = 1$ $(k = 2, .., m)$.

4.3. **Case** $m = a = 0$. In this case, the required surjective homomorphism exists if and only if

$$gcd(n/n_1, ..., n/n_b) = 1$$

by setting $\Phi(E_j) = x^{n/n_j}$.

□

Every general group which is a $\mathbb{Z}_n$-Schottky group will be called a *geometrical $\mathbb{Z}_n$-Schottky group*.

**Remark 14.** If $K$ is a general group of type $(m, a, b, n_1, ..., n_b, l_1, ..., l_m)$ and $\Omega$ is its region of discontinuity, then $\Omega/K$ is an orbifold of signature $(m + a, 2b; n_1, n_1, n_2, n_2, ..., n_b, n_b)$. In particular, if $K$ is a geometrical $\mathbb{Z}_n$-Schottky group and $\Gamma \triangleleft K$ is a Schottky group of rank $g$ so that $K/\Gamma \cong \mathbb{Z}_n$, then (by the Riemann-Hurwitz formula) we see that

$$g = n(m + a - 1) + 1 + \sum_{j=1}^{b} \frac{n}{n_j}(n_j - 1).$$

If $M = (\mathbb{H}^3 \cup \Omega)/\Gamma$ (a handlebody of genus $g$) and $\tau$ is the automorphism induced by $K$ on $M$, then the locus of fixed points of $\tau$ has as connected components exactly $\sum_{j=1}^{b} n/n_j$ simple geodesic arcs and $\sum_{k=1}^{m} n/l_k$ simple closed geodesics.

4.4. **Classification of $\mathbb{Z}_n$-Schottky groups.** In the previous section we have constructed the so called geometrical $\mathbb{Z}_n$-Schottky groups. In this section we prove every $\mathbb{Z}_n$-Schottky groups is a geometrical one; so Theorem 1 will be consequence of Theorem 13.

Let $K$ be a $\mathbb{Z}_n$-Schottky group of genus $g \geq 2$ and $\Gamma$ be a Schottky normal subgroup so that $H = \langle \tau \rangle = K/\Gamma \cong \mathbb{Z}_n$. Both groups, $K$ and $\Gamma$, have the same region of discontinuity, say $\Omega$. As $\Gamma$ has neither parabolic transformations nor elliptic transformations, the group $K$ has neither parabolic transformations nor elliptic transformations of order different from a divisor of $n$.

We consider a Schottky uniformization $(\Omega, \Gamma, P : \Omega \to S)$. We have that $H < Aut^+(S)$ lifts with respect to this Schottky uniformization to the group $K$. As a consequence of Theorem 10, we may find a fundamental collection of loops, say $\mathcal{F}$, associated to the pair $\{(\Omega, \Gamma, P : \Omega \to S), H\}$. Let us denote by $\widehat{\mathcal{F}}$ the collection of pairwise disjoint simple loops in $\Omega$ obtained by the lifting of the loops in $\mathcal{F}$ under $P$. Each of these lifted loops is called a *structure loop* for $K$ and each of the regions in $\Omega - \widehat{\mathcal{F}}$ is called a *structure region* for $K$.

There is a regular holomorphic branched covering $Q : \Omega \to O = \Omega/K$ and a regular holomorphic branched covering $\pi : S \to O$, with $H$ as deck group of covering transformations, so that $Q = \pi \circ P$.

**Lemma 15.** *The stabilizer in $K$ of either a structure loop or a structure region is either trivial or a cyclic group of order some divisor of $n$.*

*Proof.* This is a consequence of the fact that $P$ restricted to a structure loop $\alpha$ (respectively, a structure region $R$) is a homeomorphism onto the image loop $P(\alpha)$ (respectively, onto its image $P(R) \subset S$). So, there is a natural isomorphism between the $H$-stabilizer of $P(\alpha)$ (respectively, $P(R)$) and the $K$-stabilizer of $\alpha$ (respectively, $R$). □

**Lemma 16.** *Let $R$ be a structure region with non-trivial stabilizer, say generated by the non-trivial elliptic transformation $E \in K$. If one of the fixed points of $E$ belongs to $R$, then the other fixed point also belongs to $R$.*

*Proof.* As consequence of Theorem 12, both fixed points of $E$ belong to $\Omega$. Let us assume one of the fixed points of $E$ belongs to $R$ and that the other fixed point does not. It follows that there is a boundary structure loop $\gamma$ of $R$ which is invariant under $E$ (this loop separates both fixed points of $E$). Let us denote by $R'$ the other structure region bounded by $\gamma$. If $P(R) = P(R')$, then there should be some boundary structure loop $\eta$ for $R$ and some $T \in K - \{I\}$ so that $T(\eta) = \gamma$. If $\eta = \gamma$, then $T$ should be an elliptic transformation of order 2 that permutes $R$ with $R'$. In this case, we have that $\gamma$ is stabilized by a dihedral group (generated by $E$ and $T$), a contradiction by Lemma 15. If $\eta \neq \gamma$, then $T$ is a loxodromic transformation and $T^{-1}ET$ is elliptic transformation, of the same order as $E$, keeping $R$ invariant. It follows that $\langle T^{-1}ET \rangle = \langle E \rangle$. But, as $E$ does not keep invariant the loop $\eta$ and $T^{-1}ET$ does it, we get a contradiction. We have proved that $P(R) \neq P(R')$. In this way, we may delete the loop $P(\gamma)$ and its $H$-translates from $\mathcal{F}$ without destroying the properties of being a fundamental collection of loops; again a contradiction with the minimality of $\mathcal{F}$. □

**Lemma 17.** *Let $R$ be a structure region with non-trivial stabilizer, say generated by the non-trivial elliptic transformation $E \in K$. If none of the fixed points of $E$ belong to $R$, then there is a loxodromic transformation $T \in K$ and (different) structure loops $\gamma$ and $\gamma'$ in the border of $R$, each of them invariant under $E$, so that $T(\gamma) = \gamma'$.*

*Proof.* If none of the fixed points of $E$ belong to $R$, then there are two different boundary structure loops of $R$, say $\gamma$ and $\gamma'$, each one invariant under $E$. Assume there is no loxodromic transformation $T \in K$ so that $T(\gamma) = \gamma'$. Let $R'$ be the other structure region bounded by $\gamma$. If $P(R) \neq P(R')$, then we may delete from $\mathcal{F}$ the projection of $\gamma$ (and all its translates under $H$) to obtain a contradiction to the minimality of $\mathcal{F}$. If $P(R) = P(R')$, then there should be some boundary structure loop $\eta$ for $R$ and some $B \in K - \{I\}$ so that $B(\eta) = \gamma$. If $\eta = \gamma$, then $B$ is an elliptic transformation of order two and $\gamma$ will be stabilized by a dihedral group (generated by $B$ and $E$), a contradiction by Lemma 15. It follows that $\eta \neq \gamma$ and that $B$ is a loxodromic transformation. By our assumption, $\eta \neq \gamma'$. In this case, the elliptic transformations $B^{-1}EB \in K$ keeps invariant $R$, in particular, $\langle B^{-1}EB \rangle = \langle E \rangle$. But $E$ keeps invariant $\gamma$ meanwhile $B^{-1}EB$ does not, a contradiction. □

**Lemma 18.** *Let $R$ be a structure region and $K_R$ be its $K$-stabilizer.*
  (1) *If $K_R$ is trivial, then each of the boundary loops of $R$ has trivial $K_R$-stabilizer.*
  (2) *If $K_R = \langle E \rangle$ is a non-trivial cyclic group and both fixed points of $E$ belong to $R$, then each of the boundary loops of $R$ has trivial $K_R$-stabilizer.*
  (3) *If $K_R = \langle E \rangle$ is a non-trivial cyclic group and both fixed points of $E$ do not belong to $R$, then each of the boundary loops of $R$ (with the exception of the two loops separating $R$ from a fixed point of $E$) has trivial $K_R$-stabilizer.*

*Proof.* This is easily follows from the action of an elliptic transformation in the Riemann sphere. □

**Remark 19.** Let us consider two different structure regions, say $R_1$ and $R_2$, both having a common boundary loop, say $\alpha$. We assume that $\alpha$ has trivial $K_{R_1}$-stabilizer. We known, from Lemma 15, that the $K$-stabilizer of $R_j$ is either trivial or a cyclic group. Let $K_j^0 = \langle E_j \rangle$ be the $K$-stabilizer of $R_j$ (the transformation $E_j$ may be the identity). If $K_j^0$ is either trivial or both fixed points of $E_j$ belongs to $R_j$, then we set $K_j = K_j^0$ (the a trivial group or one of type (T2)). If none of the fixed points of $E_j$ belong to $R_j$, then by Lemma 17 there is a loxodromic transformation $T_j \in K$ commuting with $E_j$. In this last case we set $K_j = \langle E_j, T_j \rangle$ (a group of type (T4)). We know that the $K$-stabilizer of $\alpha$ can only be trivial or a cyclic group generated by an involution permuting the regions $R_1$ and $R_2$. Now, as a consequence of the Klein-Maskit combination theorems, we have one of the following facts.
  (i) If $\alpha$ has trivial $K$-stabilizer, then the group $\langle K_1, K_2 \rangle$ is the free product of $K_1$ with $K_2$ (applying (1) in Klein-Maskit combination theorem).
  (i) If $\alpha$ is stabilized by an involution, say $\eta \in K$, that permutes both regions, then the group $\langle K_1, K_2 \rangle$ is the HNN-extension of $K_1$ by $\langle \eta \rangle$.

4.5. **Case 1.** Assume there is a structure region $R$ with trivial $K$-stabilizer. Then, $Q : \Omega \to O$ restricted to $R$ is a homeomorphism onto its image $Q(R)$.

Let $\gamma$ be a structure loop on the boundary of $R$. If the stabilizer of $\gamma$ in $K$ is generated by an elliptic transformation of order two, say $E$, then $E$ permutes both structure regions bounded by $\gamma$. We obtain, in this way, a cyclic group of order two generated by the involution $E$.

If the stabilizer of $\gamma$ in $K$ is trivial, then there is another structure boundary loop $\gamma'$ of $R$ (different from $\gamma$) and some loxodromic transformation $B \in K$ so that $B(\gamma) = \gamma'$. In fact, if this is not the case and $R'$ denotes the other structure region bounded by $\gamma$, then $P(R) \neq P(R')$. In this case one may delete $P(\gamma)$ and $\tau^r(P(\gamma))$, $r = 1, ..., n-1$, from $\mathcal{F}$ without destroying the invariance and the planar condition, a contradiction to the minimality of $\mathcal{F}$. We have obtained a cyclic group generated by the loxodromic transformation $B$.

Working out first with all the boundary loops which have stabilizer an elliptic transformation of order two, we obtain a Kleinian group $K_1 < K$ as the free product of these cyclic groups of order two (using (1) in the Klein-Maskit combination theorem). Next, working out with the rest of the boundary loops, we obtain a free product of the previous group $K_1$ with cyclic groups generated by loxodromic transformations (again using (1) in the Klein-Maskit combination theorem), obtaining a Kleinian group $K_2 < K$. It is clear that the group $K_2$ is constructed by using groups of types (T1+) and (T2).

We may also notice from the above that every other structure region is $K_2$-equivalent to $R$ and, in particular, $K = K_2$.

4.6. **Case 2.** Assume, from now on, that every structure region has non-trivial stabilizer in $K$. We may consider a finite collection of structure regions, say $R_1, ..., R_l$, so that they are pairwise non-equivalent under the action of $K$ and maximal with such a property. It is not difficult to see that such a collection of structure regions can also be assumed to satisfies that the union of the closure of them is connected (this is just a consequence of the fact that $O$ is connected). For each region $R_j$ we consider the group $K_j < K$ as follows.

(1) If the $K$-stabilizer of $R_j$ is trivial, then set $K_j$ the trivial group.
(2) If the $K$-stabilizer of $R_j$ is a non-trivial cyclic group, generated by an elliptic transformation $E_j \in K$, for which both fixed points belong to $R_j$, then we set $K_j$ the $K$-stabilizer of $R_j$ (a group of type (T2)).
(3) Assume the $K$-stabilizer of $R_j$ is a non-trivial cyclic group, generated by an elliptic transformation $E_j \in K$, for which both fixed points do not belong to $R_j$. In this case, we know the existence of a loxodromic transformation $T_j \in K$ commuting with $E_j$ (the transformation $T_j$ send one of the boundary loops of $R_j$ that separates the fixed points onto the other boundary loop with the same property). We set $K_j = \langle E_j, T_j \rangle$ (a group of type (T4)).

Let us consider two of the regions $R_i$ and $R_j$ with a common boundary loop, say $\alpha_{ij}$. Clearly the loop $\alpha_{ij}$ cannot have non-trivial stabilizer in $K_i$ (neither $K_2$). In fact, if this happen, then the only possibility is that $\alpha_{ij}$ is a boundary loop of $R_j$ that separates the fixed points of the elliptic transformation $E_i$. But in that case, we have the loxodromic transformation $T_j$ (or $T_i^{-1}$) that sends it to the other boundary loop of $R_i$ that separates the fixed points of $E_i$; which means that $R_i$ and $R_j$ are $K$-equivalent, a contradiction. It then follows that we may perform the free product (part (1) in the Klein-Maskit combination theorem) between $K_1$ and $K_2$ (see Remark 24).

Let us denote by $R_0$ the union of all closures of the regions $R_1, ..., R_l$, and let us denote by $K_0 < K$ the Kleinian group obtained by the above procedure (working at all pairs of regions with a common boundary). By the construction, the group $K_0$ is constructed from groups of type (T2) and (T4).

Next, for every structural loop $\delta$ in the boundary of $R_0$ we have only two possibilities: either $\delta$ is invariant under an elliptic element of order two or there is another boundary loop $\delta'$ of $R_0$ and a loxodromic transformation $T \in K$ with $T(\delta) = \delta'$. In the first case we perform a free product of $K_0$ with the cyclic group generated by $T$ and in the second one we perform the free product with the cyclic group generated by $T$. Proceeding with all boundary loops of $R_0$, we ends with a Kleinian group $K^* < K$, constructed by using groups of types (T1), (T2) and (T4). As every structural region of $K$ is $K^*$-equivalent with one contained inside the region $R_0$, it follows that $K = K^*$.

## 5. Topological classification of extended $\mathbb{Z}_n$-Schottky groups

The arguments are very similar to the ones used in the proof of Theorem 1, so we only provide the main steps without going into many details. Let $K$ be a extended $\mathbb{Z}_n$-extended Schottky group of rank $g$ and $\Gamma \triangleleft K$ be a Schottky group of rank $g \geq 2$ so that $H = K/\Gamma = \langle \tau \rangle \cong \mathbb{Z}_{2n}$. The Riemann surface $S = \Omega/\Gamma$ admits the anticonformal automorphism $\tau$ of order $2n$ which lifts to the region of discontinuity $\Omega$ as an extended Möbius transformation $\widehat{\tau}$ which normalizes $\Gamma$, $\widehat{\tau}^{2n} \in \Gamma$ and $\widehat{\tau}^j \notin \Gamma$, for $j = 1, ..., 2n-1$.

We consider a Schottky uniformization $(\Omega, \Gamma, P : \Omega \to S)$. We have that $H < Aut(S)$ lifts with respect to this Schottky uniformization to $K$. As a consequence of Theorem 10, we may consider a fundamental collection $\mathcal{F}$ associated to the pair $\{(\Omega, \Gamma, P : \Omega \to S), H\}$. Let us denote by $\widetilde{\mathcal{F}}$ the collection of pairwise disjoint simple loops in $\Omega$ obtained by the lifting of the loops in $\mathcal{F}$ under $P$. Each of these lifted loops is called a *structure loop* for $K$ and each of the regions in $\Omega - \widetilde{\mathcal{F}}$ is called a *structure region* for $K$.

Similarly as in the proof of Theorem 1, we may determine the stabilizers of each structure region and a maximal collection of non-equivalent structure regions (means that two different structure regions are non-equivalent under $K$, but it may happen that the stabilizer of any of the structures regions is non-trivial in $K$).

It is possible to obtain the following general picture.

**Proposition 20.** *Every extended $\mathbb{Z}_n$-Schottky groups is obtained by the Klein-Maskit free product combination of groups of the following type:*

1.- *cyclic group generated by a glide-reflection transformation;*
2.- *cyclic group generated by a loxodromic transformation;*
3.- *cyclic group generated by elliptic transformation of order a divisor of $n$;*
4.- *cyclic group generated by a pseudo-elliptic transformation of order $2d$, where $2d$ is a divisor of $2n$, but not a divisor of $n$;*
5.- *cyclic group generated by a reflection (only if $n$ is odd);*
6.- *Abelian group generated by a loxodromic transformation and an elliptic transformation of order a divisor of $n$ (in particular, both fixed points of the elliptic are the same as for the loxodromic);*
7.- *a group generated by a loxodromic transformation $A$ and a pseudo-elliptic transformation $B$ of order a divisor of $2n$ but not of $n$ so that $B^{-1}ABA = I$ (in particular, both fixed points of the loxodromic transformation are permuted by the pseudo-elliptic transformation);*
8.- *a group generated by a glide-reflection transformation $A$ and an elliptic transformation $B$ of order a divisor of $n$ so that $B^{-1}A = AB = I$ (in particular, both fixed points of the elliptic are the same as for the glide-reflection);*
9.- *a group generated by a glide-reflection transformation $A$ and a pseudo-elliptic transformation $B$ of order a divisor of $2n$ but not of $n$ so that $(B^{-1}A)^2 = (BA)^2 = I$ (in particular, both fixed points of the glide-reflection are permuted by the pseudo-elliptic);*
10.- *groups generated by a discrete group $F$ (of orientation-preserving conformal automorphisms) keeping invariant a circle $\Sigma$ and the reflection on $\Sigma$, so that $\Omega(F)/F$ is a connected Riemann orbifold with conical points of orders which are divisors of $n$ (only if $n$ is odd).*

A group $K$ constructed with the groups listed in Proposition 20 by use of the Klein-Maskit combination theorems will be called a *general group*. As we will see below, there are general groups which are not a extended $\mathbb{Z}_n$-Schottky group.

If $K$ is a general group, then, as consequence of the Klein-Maskit combination theorems, the following properties can be observed.

(1) $K$ is a discrete group;
(2) the limit set of $K$ is a Cantor set;
(3) $K^+$ contains no parabolic transformations;
(4) every non-loxodromic transformation in $K^+$ is either the identity or conjugated to a power of some elliptic or pseudo-elliptic generator used in the basic groups construction.

The following gives a necessary and sufficient condition for a general group to be a extended $\mathbb{Z}_{2n}$-Schottky group.

**Proposition 21.** *A general group $K$ is a extended $\mathbb{Z}_n$-Schottky group if and only if*

(1) *$K$ contains orientation-reversing transformations, and*
(2) *there is a surjective homomorphism*

$$\Phi : K \to \mathbb{Z}_{2n} = \langle x : x^{2n} = 1 \rangle$$

*with torsion free kernel containing only orientation-preserving transformations.*

*Proof.* One direction is clear. The other one is consequence of the following. If $\Gamma$ is the kernel of $\Phi$ as required, then it is a function group containing only loxodromic transformations and whose limit set is a Cantor set. Then, as consequence of the classification of function groups [19] $\Gamma$ is a Schottky group. □

**Remark 22.** Condition (1) in Proposition 21 is equivalent to say that in the construction of the general group $K$, by use of the Klein-Maskit's combination theorems, we need to use at least one of the basic groups given in Proposition 20 of types 1.-, 4.-, 5.-, 7.-, 8.-, 9.- or 10.-. Once we have condition (1) full-filled, condition (2) is trivial if we have used in the construction basic groups of types 1.-, 8.- or 9.-. Let us assume we don't use any of these three types of groups. In this case, in order for to have condition (2) full-filled, we need that the maximum common divisor of all values of the form $2n/r$, were $r$ runs over all orders of elliptic and pseudo-elliptic transformations used in the basic groups, is 1.

5.1. **Elimination of groups of type 9.-.** Assume we use a group of type 9.- to construct a $\mathbb{Z}_{2n}$-extended Schottky group, say generated by a glide-reflection $A$ and a pseudo-elliptic $B$ so that the order of $B$ divides $2n$ but not $n$ and $(BA)^2 = (B^{-1}A)^2 = I$. Proposition 21 asserts the existence of a surjective homomorphism $\Phi : K \to \mathbb{Z}_4 = \langle x : x^{2n} = 1 \rangle$ with kernel $\Gamma$ which is torsion free and only containing orientation-preserving automorphisms. As $(BA)^2 = (B^{-1}A)^2 = I$ and $\Gamma$ is torsion free, we should have that $\Phi(BA) = x^n = \Phi(B^{-1}A)$, from which we obtains that $\Phi(B^2) = 1$, that is, $B^2 \in G$ and, in particular, $B^2 = I$. Now, $\Phi(B) = x^n$ as $\Gamma$ is torsion free and $B \neq I$. As $\Phi(BA) = x^n$, we should have $\Phi(A) = 1$ and $A \in \Gamma$, a contradiction to the fact that $\Gamma$ only contains orientation-preserving transformations.

5.2. **Elimination of groups of type 8.- for $n$ odd.** Assume we use a group of type 8.- to construct a $\mathbb{Z}_{2n}$-extended Schottky group, say generated by a glide-reflection $A$ and an elliptic $B$ so that the order of $B$ divides $n$ and $B^{-1}A = AB$. Proposition 21 asserts the existence of a surjective homomorphism $\Phi : K \to \mathbb{Z}_4 = \langle x : x^{2n} = 1 \rangle$ with kernel $\Gamma$ which is torsion free and only containing orientation-preserving automorphisms. As $AB = B^{-1}A$ we obtains that $\Phi(B^2) = 1$, that is, $B^2 \in \Gamma$ and, in particular, $B^2 = I$. Now, as $\Gamma$ is torsion free and $B \neq I$, necessarily $\Phi(B) = x^n$. Let $T \in K$ be so that $\Phi(T) = x$. We claim that $T \in K - K^+$ and that $n$ is even. In fact, if $T \in K^+$ and $U \in K - K^+$, then $\Phi(U) = \Phi(T^r)$ some $r \in \{1, 2, ..., 2n - 1\}$. It follows that $U^{-1}T^r \in \Gamma < K^+$, a contradiction. Now, as $\Phi(T^n) = x^n = \Phi(B)$, one has that $T^n B \in \Gamma < K^+$, so $T^n \in K^+$, that is, $n$ is necessarily even.

5.3. **Elimination of groups of type 2.-.** As already noted in Remark 22, in order to construct a $\mathbb{Z}_{2n}$-extended Schottky group we should use at least some of the groups of type 1.-, 4.- 6.- or 7.- (we know from the above that groups of type 8.- cannot be used and groups of type 7.- can only be used for $n$ even and the elliptic generator is an involution). In this way, it is not hard to see that we may change each loxodromic generator of a group of type 2.- by one glide-reflection generator for a new group of type 1.-.

All the above complete the proof of Theorem 4.

## 6. Proof of Theorem 6

Let $p$ be a prime integer and $K$ be a $\mathbb{Z}_p$-Schottky group of rank $g$ constructed as the free product (by Klein-Maskit's combination theorem) of "$a$" cyclic groups generated by loxodromic transformations $A_1,..., A_a$, "$b$" cyclic groups generated by elliptic transformations $E_1,..., E_b$, all of them of order $p$, and "$m$" Abelian groups $H_1, ..., H_m$, where $H_j \cong \mathbb{Z} \oplus \mathbb{Z}_p$ is generated by a loxodromic transformation $T_j$ and an elliptic transformation $F_j$ of order $p$, so that $T_j F_j = F_j T_j$, with the conditions that

$$g - 1 = p(m + a + b - 1) + 1 - b.$$

The above values of $a$, $b$ and $m$ are uniquely determined by $K$ (see Remark 2).

A geometrical automorphism of $K$ is an automorphism of it induced by self-conjugation of $K$ by an orientation-preserving homeomorphism of the Riemann sphere. It is not difficult to see that a geometrical automorphism $\Psi$ of $K$ can only permute the generators $E_1,..., E_b$ (up to conjugation by some element of $K$ and inversion) and can only permute the generators $F_1,..., F_m$ (up to conjugation by some element of $K$ and inversion).

In order to obtain those Schottky groups $\Gamma \triangleleft K$ of rank $g$ so that $K/\Gamma \cong \mathbb{Z}_p$, we only need to obtain those surjective homomorphisms $\Phi : K \to \mathbb{Z}_p$ with torsion free kernel (these will be the desired Schottky groups). The torsion free condition is equivalent to $\Phi(E_j)$ and $\Phi(F_k)$ belonging to the complement of the kernel. Clearly, $\Phi$ is uniquely determined by its kernel up to composition of an automorphism of $\mathbb{Z}_p$ at the left. Let $\Gamma_1, \Gamma_2$ be two Schottky groups, both of which are normal subgroups of index $p$ in $K$. Let $\Phi_j : K \to \mathbb{Z}_p$ surjective homomorphisms with $\ker(\Phi_j) = \Gamma_j$. Then $\Psi(\Gamma_1) = \Gamma_2$ for a geometrical automorphism $\Psi$ of $K$ if and only if the corresponding surjective homomorphisms satisfy the equality $A \circ \Phi_1 = \Phi_2 \circ \Psi$, for some $A \in \text{Aut}(\mathbb{Z}_p) \cong \mathbb{Z}_{p-1}$. Now, by the above observations, we only need to count how many Schottky groups are inside $K$ up to conjugation by a geometrical automorphism of it.

Let us consider a surjective homomorphism $\Phi : K \to \mathbb{Z}_p = \langle x : x^p = 1 \rangle$ whose kernel is torsion free. By the torsion free condition, we should have that $\langle \Phi(E_j) \rangle = \langle x \rangle = \langle \Phi(F_k) \rangle$, for every $j = 1, ..., b$ and $k = 1, ..., m$. We may replace $T_k$ by the new generator $F_k^{n_k} T_k$ for a suitable $n_k$ in order to assume $\Phi(T_k) = 1$, for $k = 1, ..., m$. Similarly, if $b > 0$ or $m > 0$, then we may replace $A_j$ by the new generator $E_1^{n_j} A_j$ or $F_1^{n_j} A_j$ for suitable $n_j$ in order to assume $\Phi(A_j) = 1$. If $b = m = 0$, then the surjectivity asserts there is some $r \in \{1, ..., a\}$ so that $\Phi(A_r) \neq 1$; we may replace $A_s$ ($s \neq r$) by the new generator $A_r^{u_s} A_s$ in order to assume $\Phi(A_s) = 1$. Note that each of the previous changes is produced by a geometrical automorphism of $K$. In this way, we may assume that either:

(1) if $b + m > 0$, then $\Phi(A_j) = \Phi(T_k) = 1$, for every $j = 1, ..., a$ and $k = 1, ..., m$;
(2) if $b = m = 0$, then $\Phi(A_j) = 1$, for every $j = 2, ..., a$ and $\Phi(A_1) \neq 1$.

### 6.1. The case $b + m > 0$.
As seen above, we only need to consider surjective homomorphisms $\Phi$ so that $\Phi(A_j) = \Phi(T_k) = 1$, for every $j = 1, ..., a$ and $k = 1, ..., m$. Let $\Phi(E_j) = x^{r_j}$ and $\Phi(F_k) = x^{s_k}$, where $r_j, s_k \in \{1, 2, ..., p-1\}$, for each $j = 1, ..., b$ and each $k = 1, ..., m$. The composition at the left by a geometrical automorphism of $K$ will change the tuples $(r_1, ..., r_b)$ and $(s_1, ..., s_m)$ into tuples $(\widehat{r_1}, ..., \widehat{r_b})$ and $(\widehat{s_1}, ..., \widehat{s_m})$, where $\widehat{r_j} \in \{r_{\sigma(j)}, p - r_{\sigma(j)}\}$ for some permutation $\sigma \in S_b$ and $\widehat{s_k} \in \{s_{\eta(k)}, p - s_{\eta(k)}\}$ for some permutation $\eta \in S_m$. Following as in [13] one obtains that the number of different surjective homomorphisms, up to composing at the right by a geometrical automorphism of $K$, is

$$M(b, m) = \binom{b + (p-3)/2}{(p-3)/2} \binom{m + (p-3)/2}{(p-3)/2}$$

### 6.2. The case $b = m = 0$.
We only have one possibility given by the surjective homomorphism $\Phi$ so that $\Phi(A_j) = 1$, for every $j = 2, ..., a$, and $\Phi_1(A_1) = x$.

All of the above permits to obtain Theorem 6.

## 7. Example: Extended $\mathbb{Z}_2$-Schottky groups

In this section we obtain the number of topologically non-equivalent extended $\mathbb{Z}_2$-Schottky groups. First, in this case, Theorem 4 has the following form.

**Corollary 23.**
- I.- *Every extended $\mathbb{Z}_2$-Schottky groups of rank g is obtained by the Klein-Maskit free product combination of groups of the following type:*
- (T1).- *cyclic group generated by a glide-reflection transformation;*
- (T2).- *cyclic group generated by an elliptic transformation of order 2;*
- (T3).- *cyclic group generated by a pseudo-elliptic transformation of order 4;*
- (T4).- *Abelian group generated by a loxodromic transformation and an elliptic transformation of order 2;*
- (T5).- *a group generated by a loxodromic transformation A and a pseudo-elliptic transformation B of order 4 so that $B^{-1}AB = I$;*
- (T6).- *a group generated by a glide-reflection transformation A and an elliptic transformation B of order 2 so that $BA = AB = I$.*
- II.- *A general group K constructed using "$a_j$" groups of type (Tj) (where $j = 1, ..., 6$) is a extended $\mathbb{Z}_2$-Schottky group of rank g if and only if the following two conditions are satisfied:*
  - 1.- $a_1 + a_3 + a_5 + a_6 > 0$; and
  - 2.- $g = 4a_1 + 2a_2 + 3a_3 + 4a_4 + 4a_5 + 4a_6 - 3$.

  *In this case, we say that $(a_1, a_2, a_3, a_4, a_5, a_6)$ is the signature of the general group K.*

It is easy to see that different signatures of $\mathbb{Z}_2$-Schottky groups of the same rank $g$ produce topologically non-equivalent extended $\mathbb{Z}_2$-Schottky groups. In this way, in order to count the number of different topological extended $\mathbb{Z}_2$-Schottky groups of rank $g$, we only need to count the number of different tuples $(a_1, a_2, a_3, a_4, a_5, a_6)$ of non-negative integers satisfying

(1) $a_1 + a_3 + a_5 + a_6 > 0$; and
(2) $g = 4a_1 + 2a_2 + 3a_3 + 4a_4 + 4a_5 + 4a_6 - 3$.

Define $n : \mathbb{N}_0 \to \mathbb{N}_0$ where $n(a)$ denotes the number of different triples $(x, y, z) \in \mathbb{N}_0^3$ so that $a = x + y + z$. Denote by $N_g$ the collection of tuples $(a, b, c, d) \in \mathbb{N}_0^4$ so that (i) $g + 3 = 4a + 2b + 3c + 4d$ and (ii) if $a = 0$, then $c \geq 1$. Then, the number of topologically different extended $\mathbb{Z}_2$-Schottky groups of rank $g$ is given by

$$\sum_{(a,b,c,d) \in N_g} n(d)$$

### 7.1. Extended $\mathbb{Z}_2$-Schotky groups of rank 1.

If we take $g = 1$, then we get the following two possibilities

$$(a_1, a_2, a_3, a_4, a_5, a_6) \in \{(1, 0, 0, 0, 0, 0), (0, 0, 0, 0, 0, 1)\}.$$

The tuple $(1, 0, 0, 0, 0, 0)$ produces an $\mathbb{Z}_4$-extended Schottky group $K_1$ which, up to conjugation by a suitable Möbius transformation, is generated by $A(z) = \lambda \overline{z}$ where $\lambda > 1$. In this case, there is exactly one Schottky group $\Gamma_1$ in $K_1$ with $K_1/\Gamma_1 \cong \mathbb{Z}_4$; which is generated by $A^4(z) = \lambda^4 z$. The group $K_1$ induces on the handlebody $M_1 = \mathbb{H}^3/\Gamma_1$ an orientation-reversing isometry of order 4 acting freely. The tuple $(0, 0, 0, 0, 0, 1)$ produces a extended $\mathbb{Z}_2$-Schottky group $K_2$ which, up to conjugation by a suitable Möbius transformation, is generated by $A(z) = \lambda \overline{z}$ and $B(z) = -z$, where $\lambda > 1$. In this case, there is exactly one Schottky group $\Gamma_2$ in $K_2$ with $K_2/\Gamma_2 \cong \mathbb{Z}_4$; which is generated by $BA^2(z) = -\lambda^2 z$. The group $K_2$ induces on the handlebody $M_2 = \mathbb{H}^3/\Gamma_2$ an orientation-reversing isometry of order 4 whose locus of fixed points is a simple closed geodesic. This also ensures that there is no handlebody of genus 1 admitting simultaneously both kind of orientation-reversing isometries of order 4.

### 7.2. Extended $\mathbb{Z}_2$-Schotky groups of rank 2.

If $g = 2$, then the only tuple is $(a_1, a_2, a_3, a_4, a_5, a_6) = (0, 1, 1, 0, 0, 0)$. The produced extended $\mathbb{Z}_2$-Schottky group $K$ is, up to conjugation by a suitable Möbius transformation, generated by $A(z) = -1/\overline{z}$ and an elliptic transformation of order 2, say $B$. In this case, there is exactly one Schottky group $\Gamma$ in $K$ with $K/\Gamma \cong \mathbb{Z}_4$; which is generated by $A^2 B$ and $A^{-1} B A^{-1}$. The group $K$ induces on the handlebody $M = \mathbb{H}^3/\Gamma$ an orientation-reversing isometry of order 4 on which it has exactly one simple closed geodesic as locus of fixed points and whose square is the hyperelliptic involution (then its locus of fixed points consists of exactly 3 pairwise disjoint simple geodesic arcs. The closed Riemann surface $S$ uniformized by $\Gamma$ (that is, the conformal boundary of $M$) corresponds to an algebraic curve of the form

$$y^2 = x(x^2 - 1)(x - b)(x + \overline{b}^{-1}),$$

where $b^4 \neq 0, 1$. The orientation-reversing automorphism of order 4 is given by

$$\tau : \begin{cases} x & \mapsto & \dfrac{-1}{\overline{x}} \\ y & \mapsto & \left(\dfrac{-b}{\overline{b}}\right)^{1/2} \dfrac{\overline{y}}{\overline{z}^3} \end{cases}$$

8. Connection with handlebodies

Let $H_1$ and $H_2$ be (finite) groups of homeomorphisms of a fixed handlebody $M$. We say that they are (weakly) *topologically equivalent* if there is an orientation-preserving self-homeomorphism of $M$ that conjugates $H_1$ onto $H_2$. In [2] Bartoszyńska provided the topological classification of involutions of handlebodies of genus two, in [13] Kalliongis-Miller characterized orientation-preserving finite group actions on handlebodies of genus $g \geq 2$ and in [12] Kalliongis-McCullough considered a topological picture of orientation reversing involutions. In all of these papers, the used method is combinatorial and 3-dimensional in nature. In [13] an explicit formulae to obtain the number of equivalence classes for cyclic groups of prime order $p$ was provided. In [3] it was studied the case of free fixed point orientation-reversing group actions on handlebodies and a classification theorem was obtained in terms of algebraic invariants that involve Nielsen equivalence.

Let $M$ be a handlebody of genus $g$. If $\Gamma$ is a Schottky group of rank $g$, then $M_\Gamma = (\mathbb{H}^3 \cup \Omega)/\Gamma$, where $\Omega$ is the region of discontinuity of $\Gamma$ is homeomorphic to $M$. We say that $\Gamma$ induces a *Schottky structure* on $M$ ($\Gamma$ induces a complete hyperbolic structure on the interior of $M$, whose injectivity radius is bounded away from zero, and also it provides a Riemann surface structure on the topological boundary of $M$). A *conformal automorphism* (respectively, *anticonformal automorphism*) of $M_\Gamma$ is an orientation-preserving (respectively, orientation-reversing) homeomorphism $f : M_\Gamma \to M_\Gamma$ whose restriction to the interior hyperbolic 3-manifold $M_\Gamma^0 = \mathbb{H}^3/\Gamma$ is an isometry.

The following well known fact permits to see the relation between finite groups of homeomorphisms of handlebodies and Schottky groups.

**Lemma 24.** *Let $M$ be a handlebody of genus $g$ and let $H$ be a finite group of homeomorphism of $M$. Then there is a Schottky structure on $M$ and there is a group $\widehat{H} \cong H$ of conformal/anticonformal automorphisms of $M$ (with respect to the Schottky structure) which is homotopic to $H$.*

*Proof.* Let $M^0$ be the interior of $M$ and let $S$ is boundary (a closed orientable surface of genus $g$). The group $H$ acts on $S$ as a group $H_S$ of homeomorphisms. As a consequence of Nielsen's realization theorem [14], we may assume that $S$ has a Riemann surface structure so that (up to homotopy equivalence) the group $H_S$ acts as a group of automorphisms on $S$. The Riemann surface structure on $S$ induces a complete hyperbolic structure on $M^0$ so that its conformal boundary is exactly the Riemann surface $S$, that is a Schottky structure on $M$. Let us assume that such a Schotty structure is provided by the Schottky group $\Gamma$. The Schottky uniformization $(\Omega, \Gamma, P : \Omega \to S)$ has the property that the group $H_S < \text{Aut}(S)$ lifts to a group $K$ of automorphisms of $\Omega$. We already know that $K$ is a subgroup of $\widehat{\mathbb{M}}$ with $\Gamma \triangleleft K$ and $H_S \cong K/\Gamma$. The group $K$ induces a group $\widehat{H}$ of conformal/anticonformal automorphisms of $M$, isomorphic to $H$, which is homotopic to $H$ on the boundary $S$. As $M$ is a compression body, $\widehat{H}$ and $H$ are homotopic on $M$. □

8.1. **The cyclic case.** Let us assume that $M$ is a handlebody of genus $g \geq 2$ and that $\tau : M \to M$ is a finite order homeomorphism and set $H = \langle \tau \rangle$. Let the order of $\tau$ be $n$ (respectively, $2n$) if $\tau$ is orientation-preserving (respectively, orientation-reversing). By Lemma 24, there is a Schottky structure on $M$ for which we may though of $\tau$ (up to homotopy) as a conformal/anticonformal automorphism. Let us denote by $\Gamma$ the Schottky group which provides such Schottky structure on $M$. Theorems 1 and 4 complement the work done by Kalliongis-Miller in [13] as follows. By lifting $\tau$ to the universal cover space (the hyperbolic 3-dimensional space), we obtain a $\mathbb{Z}_n$-Schottky group (respectively, extended $\mathbb{Z}_n$-Schottky group) $K$ of rank $g$ if $\tau$ is conformal (respectively, anticonformal) with $\Gamma \triangleleft K$ and $K/\Gamma = \langle \tau \rangle$.

8.1.1. *Case $\tau$ is conformal.* In this case, as $K$ is a $\mathbb{Z}_n$-Schottky group of rank $g$, the structural decomposition provided by Theorem 1 asserts that $K$ can be constructed using "$a$" cyclic groups, each one generated by a loxodromic transformation, "$b$" cyclic groups generated by elliptic transformations and "$m$" Abelian groups. Assume the $b$ elliptic cyclic groups have orders $n_1,..., n_b$, and that the $m$ Abelian groups are isomorphic to $\mathbb{Z}_{l_1} \oplus \mathbb{Z},..., \mathbb{Z}_{l_m} \oplus \mathbb{Z}$. With this information we are able to describe the locus of fixed points as $\tau$ and the quotient orbifold $M/H$ as follows.

(1) The locus of fixed points of $\tau$ is given by a pairwise collection of $\sum_{j=1}^{m} n/l_j$ simple loops (closed geodesics of $M_\Gamma^0$) and $\sum_{j=1}^{b} n/n_j$ simple arcs connecting two different points on $S$ (the interiors of these arcs being simple geodesics of $M_\Gamma^0$). For instance, as can be seen from Section 2.4.1, an automorphism of order $n = 11$ acting on a handlebody of genus $g = 100$ has no loops of fixed points, by the exception of one case. In the exceptional case the involution has exactly 11 components of fixed points, each one being a loop.
(2) The quotient orbifold $M/H$ is a (topological) handlebody of genus $b+m$ whose conical locus consists of exactly "$b$" simple arcs and "$m$" simple loops; all of them disjoint.

If $n = 2$, then in [13] it was noted that the number of topologically non-equivalent actions is $(g+2)(g+4)/8$ when $g$ is even and $(g + 3)(g + 5)/8$ when $g$ is odd. This number is the same as the number of topological different conjugacy classes of $\mathbb{Z}_2$-Schottky groups of a fixed rank $g$. This is consequence of the fact that given any index two Schottky subgroups, say $\Gamma_1$ and $\Gamma_2$, of the same $\mathbb{Z}_2$-Schottky group $K$, then one may construct an orientation-preserving homeomorphism (even a quasi-conformal one) $f : \widehat{\mathbb{C}} \to \widehat{\mathbb{C}}$ so that $fKf^{-1} = K$ and $f\Gamma_1 f^{-1} = \Gamma_2$.

As can be seen from Section 2.4.1, the number of topologically non-equivalent $\mathbb{Z}_{13}$-Schottky groups of rank 157 is 16; but in [13] it was obtained that there are 87,108 topologically non-equivalent actions of $\mathbb{Z}_{13}$ as group of orientation-preserving self-homeomorphisms of handlebodies of genus 157. This difference corresponds to the fact that a $\mathbb{Z}_{13}$-Schottky group $K$ may contains many different Schottky groups as normal subgroups of index 13, up to a geometric automorphisms of $K$. With respect to this problem, Theorem 6 provides the number of Schottky groups of index $p$ ($p$ a prime) inside a $\mathbb{Z}_p$-Schottky group $K$ which are non-equivalent under a geometric automorphism of $K$. In this way, the counting process done in the Section 2.4 permits to re-obtain the counting formula obtained in [13] of topologically non-equivalent (up to homotopy) cyclic groups of orientation-preserving homeomorphisms of order a prime $p$ on handlebodies.

8.1.2. *Case $\tau$ is anticonformal.* In this case $K$ is a extended $\mathbb{Z}_{2n}$-Schottky group of rank $g$, so Theorem 4 asserts that $K$ can be constructed using groups of types (T1)-T(8). The description of the connected components of fixed points of $H$ is as follows.

(1) A group of type (T2), say generated by an elliptic transformation of order $d \geq 2$ a divisor of $n$, produces a simple arc as a conical component of $M/H$ (with the exception of its ends points, its lies in the interior of $M/H$). That arc has conical order equal to $d$. By lifting to $M$ such a conical arc, we obtain a collection of $2n/d$ simple arcs in $M$ (with the exception of its ends points, they lie in the interior). Each of these arcs is a connected component of the locus of fixed points of $\tau^{2n/d}$.
(2) A group of type (T3), say generated by a pseudo-elliptic transformation of order 2 (that is, an imaginary reflection) produces an isolated conical point in the interior of $M/H$ of conical order 2. It liftings to $M$ provides a collection of $n$ points, being fixed points of $\tau^n$ (in this case, $n$ is necessarily odd).
(3) A group of type (T3), say generated by a pseudo-elliptic transformation of order $2d$, with $d \geq 2$ a divisor of $n$, produces a simple arc as a conical component of $M/H$. One of the end points belong to the interior of $M/H$ and the other to the boundary. The arc (with the interior end point deleted) has conical order $d$. The deleted point has conical order $2d$. By lifting such an arc, we obtain a collection of $n/d$ simple arcs as components of fixed points of $\tau^{2n/d}$. The lifting of the interior end point is a collection of $n/d$ fixed points of $\tau^{n/d}$ (in particular, this case only happens if $n/d$ is odd).

(4) A group of type (T4), say with the elliptic generator being of order $d \geq 2$ a divisor of $n$, produces a simple loop as a conical component of $M/H$. Its conical order is $d$. By lifting such a loop, we obtain a collection of $2n/d$ simple loops as components of fixed points of $\tau^{2n/d}$.

(5) A group of type (T5), say with the pseudo-elliptic generator of order 2, produces two isolated conical point in the interior of $M/H$ of conical order 2. Its liftings to $M$ provides a collection of $2n$ points, being fixed points of $\tau^n$ (in this case, $n$ is necessarily odd).

(6) A group of type (T5), say that the pseudo-elliptic generator has order $2d$, with $d \geq 2$ a divisor of $n$, produces a simple arc as a conical component of $M/H$ of conical order $d$ (both of its end points has conical order $2d$). By lifting such an arc to $M$ we obtain a collection of $2n/d$ simple loops as components of fixed points of $\tau^{2n/d}$ and the lifted points (of the conical point of order $2d$) are isolated fixed points of $\tau^{n/d}$ (in which case, $n/d$ should be odd).

(7) A group of type (T6) (then $n$ is even) produces a simple loop as a component of conical points of $M/H$, with conical order 2. By lifting such a loop we obtain a collection of $n$ simple loops in $M$ as components of fixed points of the conformal involution $\tau^n$.

(8) A group of type (T7) (then $n$ is odd) produces a closed disc as a conical components of $M/H$ (the boundary of the disc belongs to the boundary of $M/H$ and the interior of the disc to the interior of $M/H$) of conical order 2. By lifting such a disc, we obtain a collection of $n$ discs as components of fixed points of the anticonformal involution $\tau^n$.

(9) A group of type (T8) (then $n$ is odd) produces a bordered compact orbifold as a component of conical points of $M/H$. By lifting it, we obtain a collection of bordered compact surfaces as components of fixed points of the anticonformal involution $\tau^n$.

It follows that, in the case that $n$ is odd, $\tau$ has no dimension two real locus of fixed points if and only if there are no groups of types (T7) nor (T8) in the construction of $K$. In this case, the anticonformal involution $\tau^n$ only has isolated fixed points at most. We discuss this case in the example below.

8.1.3. *Example.* Let us assume $n$ is odd and there are no groups of type (T7) nor (T8) in the construction of $K$. As $n$ is odd, then neither we cannot use groups of type (T6). Assume that in the construction of $K$ we use

(1) $a_1$ groups of type (T1);
(2) $a_2$ groups of type (T2), of orders $2 \leq l_1 \leq l_2 \leq \cdots \leq l_{a_2} \leq n$, where each $l_j$ is a divisor of $n$;
(3) $a_3$ groups of type (T3), or orders $2 \leq 2r_1 \leq 2r_2 \leq \cdots \leq 2r_{a_3} \leq 2n$, where each $2r_j$ is a divisor of $2n$ but not a divisor of $n$;
(4) $a_4$ groups of type (T4); and
(5) $a_5$ groups of type (T5).

It follows that $O = \Omega/K$ is a closed Klein orbifold of signature

$$(2a_1 + a_3 + 2a_4 + 2a_5; -; l_1, l_1, l_2, l_2, ..., l_{a_2}, l_{a_2}, r_1, r_2, ..., r_{a_3}),$$

that is, $O$ is the connected sum of $2a_1 + a_3 + 2a_4 + 2a_5$ projective planes and $2a_2 + a_3$ conical points (whose orders are given in the third part of the signature). In this way, $O^+ = \Omega/K^+$ is a Riemann orbifold of signature

$$(2a_1 + a_3 + 2a_4 + 2a_5 - 1; +; l_1, l_1, l_1, l_1, ..., l_{a_2}, l_{a_2}, l_{a_2}, l_{a_2}, r_1, r_1, ..., r_{a_3}, r_{a_3}),$$

that is, $O^+$ is a closed Riemann surface of genus $2a_1 + a_3 + 2a_4 + 2a_5 - 1$ and $4a_2 + 2a_3$ conical points.

Let $\Gamma$ be a Schottky group so that $\Gamma \triangleleft K$ and $K/\Gamma \cong \mathbb{Z}_{2n}$. The closed Riemann surface $S = \Omega/\Gamma$ admits an anticonformal automorphism $\tau|_S$, of order $2n$, induced by $\tau$. Set $f = \tau^2$. Then, $f : S \to S$ is a conformal automorphism of order $n$, with $S/\langle f \rangle = O^+$. It follows from Riemann-Hurwitz formula that the genus of $S$, that is, the rank of $\Gamma$ is

$$g = n\left(2a_1 + a_3 + 2a_4 + 2a_5 - 2 + 2\sum_{j=1}^{a_2}(1 - l_j^{-1}) + \sum_{l=1}^{a_3}(1 - r_l^{-1})\right) + 1.$$

Departamento de Matemática, Universidad Técnica Federico Santa María, Valparaiso–Chile
*E-mail address*: ruben.hidalgo@usm.cl